\documentclass[12pt,eqright]{amsart}
\usepackage[displaymath,pagewise]{lineno}
\usepackage{amssymb,amsmath}
\usepackage{bm}
\usepackage{amsfonts}
\usepackage{epsf}
\usepackage{graphicx}
\usepackage{placeins}
\usepackage{booktabs}
\newtheorem{theorem}{Theorem}[section]
\newtheorem{lemma}[theorem]{Lemma}

\newtheorem{remark}[theorem]{Remark}

\theoremstyle{definition}

\numberwithin{equation}{section}

\newcommand{\bc}{\begin{center}}
\newcommand{\ec}{\end{center}}
\newcommand{\be}{\begin{eqnarray}}
\newcommand{\ee}{\end{eqnarray}}
\newcommand{\nn}{\nonumber}
\newcommand{\ben}{\begin{eqnarray*}}
\newcommand{\een}{\end{eqnarray*}}

\newcommand{\Om}{\Omega}

\newcommand{\pa}{\partial}

\newcommand{\na}{\nabla}

\def\na{\nabla}

\def\cT{\mathcal{T}}
\def\R{\mathbb{R}}
\DeclareMathOperator{\sspan}{span}

\begin{document}
\title[]
{A lower bound of the  $L^2$ norm error estimate  for the Adini element of the
biharmonic equation}
\author[J. Hu]{Jun Hu}
\address{LMAM and School of Mathematical Sciences,
 Peking University, Beijing 100871, P. R. China}
\email{hujun@math.pku.edu.cn}

\author[Z. C.  Shi]{Zhongci Shi}
\address{LSEC, ICMSEC, Academy of Mathematics and Systems
Science, Chinese Academy of Sciences, Beijing 100190, China.}
\email{shi@lsec.cc.ac.cn}
\date{\today}
\thanks{The research of the first author was supported
 by the NSFC Project 11271035, and  in part by the NSFC Key Project 11031006.}

\maketitle
\begin{abstract}
This paper is devoted to the $L^2$ norm error estimate of the
Adini element  for the biharmonic equation. Surprisingly,  a lower bound  is
established which proves  that the $ L^2$ norm convergence rate  
can not be higher than that in  the energy norm. 
This proves the conjecture of [Lascaux and  Lesaint,
Some nonconforming finite elements for the plate bending problem,
RAIRO Anal. Numer. 9 (1975), pp. 9--53.] that the convergence rates in both $L^2$ and $H^1$ norms
can not be  higher than that in the energy norm  for this element.
\end{abstract}
\section{Introduction}
For the numerical analysis of  finite element methods for 
fourth order elliptic problems, one unsolved fundamental problem  is
the  $L^2$ norm error estimates \cite{BrennerScott,CiaBook,ShiWang10}.  In a recent paper \cite{HuShi12},
 we analyzed several mostly popular lower order elements:
 the  Powell-Sabin $C^1-P_2$ macro element \cite{PowellSabin1977}, the nonconforming
Morley element \cite{Morley68,ShiWang10,WX06,WangXu12}, the $C^1-Q_2$ macro element \cite{HuHuangZhang2011},
the nonconforming rectangle Morley element \cite{WangShiXu07}, and  the nonconforming incomplete
biquadratic element \cite{Shi86,WuMaoqing1983}. In particular, we proved that the 
best $L^2$ norm error estimates  for these elements were at most of second order
and could not be two order higher than that in the energy norm.

The Adini element \cite{BrennerScott,CiaBook,Lascaux85,ShiWang10} is one of the earliest  finite elements, dating
back over 50 years.  It is a nonconforming finite element for the biharmonic equation on rectangular meshes.
The shape function space contains the complete cubic space and two additional monomials on each rectangle.
In 1975, Lascaux and Lesaint analyzed this element and showed that the consistency error
 was of second order for uniform meshes, the same as the approximation error, and thus obtained a second order convergence rate
\cite{Lascaux85}; see also \cite{LuoLin04,MaoChen2005} and \cite{HuHuang11,Yang2000}.
 This in particular  implies at least a second order $H^1$ norm convergence rate.
 Lascaux and Lesaint also conjectured that it did not seem possible to improve this estimate; see \cite[Remark 4.5]{Lascaux85}. However, they did not provide a rigorous proof or justification  for this remark.

The purpose of this paper is to analyze the $L^2$ norm  error estimate for the Adini element \cite{BrennerScott,CiaBook,Lascaux85,ShiWang10}.
There are two main ingredients for the analysis.  One is  a refined property of the canonical interpolation operator,
which  is  proved  by  a  new  expansion method.  The other is  an identity for $(-f, e)_{L^2(\Omega)}$, where $f$ is the right-hand side function and  $e$ is the error.  Such an identity separates the dominant term from the other higher order terms, which is the key to
   use the aforementioned refined  property of the interpolation operator.
   Based on these factors,  a lower bound of the $L^2$ norm error estimate  is surprisingly
established which proves  that the best $L^2$ norm error estimate 
is at most of order $\mathcal{O}(h^2)$. Thus,  by the usual Poincare inequality, this
indicates that the best $H^1$ norm error estimate  is also  at most of order $\mathcal{O}(h^2)$.
This gives a rigorous proof of the conjecture from \cite{Lascaux85} that the convergence rates in both $L^2$ and $H^1$ norms
can not be higher than that in the energy norm.

 The paper is organized as follows.  In the following section, we
  present  the   Adini  element and define the canonical interpolation operator. 
  In Section 3,  based on  a refined property of the canonical interpolation operator and 
  an identity for $(-f, e)_{L^2(\Omega)}$,   we prove the main result that the $L^2$ norm error estimate  has a lower bound which indicates
 that the  convergence rates in both $L^2$ and $H^1$ norms  are at most of order $\mathcal{O}(h^2)$.
  In Section  4, we  analyze the  refined property of the canonical interpolation operator.  In Section 5, 
  we   establish  the identity for $(-f, e)_{L^2(\Omega)}$.  In Section 6, we end this paper by the conclusion and some comments.

\section{The  Adini element method}
We consider the model fourth order elliptic problem: Given $f\in
L^2(\Omega)$ find $w\in W:=H_0^2(\Om)$, such that
\begin{equation}\label{cont}
\begin{split}
a(w, v):=(\na^2 w, \na^2 v)_{L^2(\Om)} =(f,v)_{L^2(\Om)}\text{ for
any } v\in W.
\end{split}
\end{equation}
where $\nabla^2 w$ denotes the Hessian matrix of the function $w$.

To consider the discretization of \eqref{cont} by the Adini element,  let $\mathcal{T}_h$ be a  uniform regular   rectangular
triangulation of the domain $\Om\subset\R^2$ in the two dimensions.
Given $K\in \cT_h$, let $(x_c,y_c)$ be the center of $K$, the
horizontal length $2h_{x,K}$,  the vertical length
 $2h_{y,K}$, which define the meshsize $h:=\max\limits_{K\in\cT_h}\max(h_{x,K},h_{y,K})$ and
 affine mapping:
 \begin{equation}\label{mapping}
 \xi:=\frac{x-x_c}{h_{x, K}},\quad \eta:=\frac{y-y_c}{h_{y, K}} \text{ for any
 } (x,y)\in K.
 \end{equation}
 Since $\cT_h$ is a uniform mesh, we define $h_x:=h_{x,K}$ and $h_y:=h_{y,K}$ for any $K$. 
  On element $K$, the shape function space of the Adini element reads
\cite{BrennerScott,CiaBook,Lascaux85,ShiWang10}
 \begin{equation}\label{nonconformingnDAd} Q_{Ad}(K):
 =P_3(K)+\sspan\{ x^3y, y^3x\}\,,
  \end{equation}
  here and throughout this paper, $P_\ell(K)$ denotes the space of
  polynomials of degree $\leq \ell$ over $K$.
   The  Adini  element space $W_h$ is then defined
by
\begin{equation}\label{AdininD}
 W_h:=\begin{array}[t]{l}\big\{v\in L^2(\Om):
 v|_{K} \in Q_{Ad}(K) \text{ for each }K\in \mathcal{T}_h,
 v \text{ and } \na v \text{ is continuous }\\[0.5ex]
 \text{ at the internal nodes,  and vanishes at the boundary nodes on
 }\pa\Om
 \big\}\,.
 \end{array}\nn
\end{equation}

The finite  element approximation of Problem \eqref{cont} reads:
Find $w_h\in W_h$, such that
\begin{equation}\label{disc}
\begin{split}
a_h(w_h,  v_h):=(\na_h^2 w_h, \na_h^2 v_h)_{L^2(\Om)}
=(f,v_h)_{L^2(\Om)}\text{ for any } v_h\in W_h\,,
\end{split}
\end{equation}
where the operator $\nabla_h^2$ is the discrete counterpart of
$\nabla^2$, which is defined element by element since the discrete
space $W_h$ is nonconforming.

Given $K\in\cT_h$,  define the  canonical interpolation operator
$\Pi_K : H^3(K)\rightarrow Q_{Ad}(K)$ by, for any $v\in H^3(K)$,
\begin{equation}
(\Pi_K v)(P)=v(P)\text{ and }\na(\Pi_K v)(P)=\na
v(P)
\end{equation}
 for any vertex $P$ of $K$.  The interpolation $\Pi_K$  has the following estimates \cite{BrennerScott,CiaBook,Lascaux85,ShiWang10}:
 \begin{equation}\label{AdinterpolationEstimate}
 |v-\Pi _Kv|_{H^{\ell}(K)}\leq C h^{4-\ell}|v|_{H^4(K)}, \ell=1,2,3,4\,,
 \end{equation}
 provided that $v\in H^4(K)$.   Herein and throughout this paper, $C$ denotes a generic positive  constant which is independent of the meshsize and
  may be different at different places.  Then the global version $\Pi_h $ of the interpolation
operator $\Pi_K $ is defined as
\begin{equation}\label{interpolation}
\Pi_h |_K=\Pi_K  \text{ for any } K\in \cT_h.
\end{equation}

\section{A lower bound of the  $L^2$ norm error estimate}

This section proves the main result of this paper, namely, a lower bound of the   $L^2$ norm error estimate.
The main ingredients are a lower bound of $a_h(w-\Pi_hw,\Pi_hw)$ in Lemma \ref{Lemma2.2} 
and an identity for $(-f, w-w_h)_{L^2(\Omega)}$ in Lemma \ref{Lemma3.1}.

For the analysis, we list two results  from \cite{Lascaux85} and \cite{LuoLin04,MaoChen2005}.

\begin{lemma}\label{Lemma3.2} Let $w\in H_0^2(\Omega)\cap H^4(\Omega)$ be the solution of problem \eqref{cont}.  It
holds that
\begin{equation}
|a_h(w,v_h)-(f,v_h)|\leq Ch^2|w|_{H^4(\Omega)}\|\nabla_h^2
v_h\|_{L^2(\Omega)}\text{ for any }v_h\in W_h.
\end{equation}
\end{lemma}

\begin{lemma}\label{Lemma3.3} Let $w$ and $w_h$ be solutions of problems \eqref{cont}
and \eqref{disc}, respectively. Suppose that $w\in H_0^2(\Omega)\cap
H^4(\Omega)$. Then,
\begin{equation}
\|\nabla_h^2(w-w_h)\|_{L^2(\Omega)}\leq Ch^2|w|_{H^4(\Omega)}.
\end{equation}
\end{lemma}

\begin{theorem}\label{main} Let $w\in H_0^2(\Omega)\cap H^4(\Omega)$ and $w_h$ be solutions of problems \eqref{cont}
and \eqref{disc}, respectively.  Then, there exists a positive constant $\alpha$
which is independent of the meshsize such that
\begin{equation}
\alpha h^2\leq \|w-w_h\|_{L^2(\Omega)},
\end{equation}
provided that $\|f\|_{L^2(\Omega)}\not=0$ and that the meshsize is small enough.
\end{theorem}
\begin{proof}  The main ingredients for the proof  are Lemma \ref{Lemma3.1} and Lemma \ref{Lemma2.2}.
 Indeed, it follows from Lemma \ref{Lemma3.1} that
 \begin{equation}\label{eq3.0a}
\begin{split}
&(-f, w-w_h)_{L^2(\Omega)}\\
&=a_h(w,\Pi_hw-w_h)-(f, \Pi_hw-w_h)_{L^2(\Omega)}\\
&\quad +a_h(w-\Pi_hw,w-\Pi_hw)+a_h(w-\Pi_hw,w_h-\Pi_hw)\\
&\quad +2(f, \Pi_hw-w)_{L^2(\Omega)}+2a_h(w-\Pi_hw,\Pi_hw).
\end{split}
\end{equation}
The first two terms on the right-hand side of \eqref{eq3.0a} can be
bounded by Lemmas \ref{Lemma3.2}-\ref{Lemma3.3},  and the estimates of
\eqref{AdinterpolationEstimate}, which leads to
\begin{equation*}
\begin{split}
&|a_h(w,\Pi_hw-w_h)-(f, \Pi_hw-w_h)_{L^2(\Omega)}|\\
& \leq
Ch^2|w|_{H^4(\Omega)}\|\nabla_h^2(\Pi_hw-w_h)\|_{L^2(\Omega)}\\
&\leq Ch^2 \big(
\|\nabla_h^2(\Pi_hw-w)\|_{L^2(\Omega)}+\|\nabla_h^2(w-w_h)\|_{L^2(\Omega)}\big)\\
&\leq Ch^4 |w|_{H^4(\Omega)}^2.
\end{split}
\end{equation*}
The estimates of the third and fifth terms on the right-hand side of
\eqref{eq3.0a} follow immediately from
\eqref{AdinterpolationEstimate}, which gives
\begin{equation*}
\begin{split}
&|a_h(w-\Pi_hw,w-\Pi_hw)+2(f, \Pi_hw-w)_{L^2(\Omega)}|\\
& \leq
Ch^4\big(|w|_{H^4(\Omega)}+\|f\|_{L^2(\Omega)}\big)|w|_{H^4(\Omega)}.
\end{split}
\end{equation*}
From the Cauchy-Schwarz inequality,  the triangle inequality, Lemma
\ref{Lemma3.3}, and \eqref{AdinterpolationEstimate} it follows that
\begin{equation*}
|a_h(w-\Pi_hw,w_h-\Pi_hw)|\leq Ch^4|w|_{H^4(\Omega)}^2.
\end{equation*}
The last term on the right hand-side of \eqref{eq3.0a} has already
been analyzed in Lemma \ref{Lemma2.2}, which reads
\begin{equation*}
\beta h^2 \leq (\nabla_h^2(w-\Pi_hw),
\nabla_h^2\Pi_hw)_{L^2(\Omega)},
\end{equation*}
for some positive constant $\beta$.  A combination of these
estimates  states
\begin{equation*}
\delta h^2\leq (-f, w-w_h)_{L^2(\Omega)},
\end{equation*}
for some positive constant $\delta$ which is independent of the
meshsize provided that the meshsize is small enough. This plus the
definition of the $L^2$ norm of $w-w_h$ prove
\begin{equation*}
\begin{split}
\|w-w_h\|_{L^2(\Omega)}&=\sup\limits_{0\not=d\in
L^2(\Omega)}\frac{(d, w-w_h)_{L^2(\Omega)}}{\|d\|_{L^2(\Omega)}}\\
&\geq \frac{(-f,
w-w_h)}{\|-f\|_{L^2(\Omega)}}\geq \delta/\|f\|_{L^2(\Omega)}h^2.
\end{split}
\end{equation*}
Setting $\alpha=\delta/\|f\|_{L^2(\Omega)}$ completes the proof.
\end{proof}

\begin{remark} By the Poincare inequality, it follows that
\begin{equation*}
\alpha h^2 \leq \|\nabla(w-w_h)\|_{L^2(\Omega)}.
\end{equation*}
\end{remark}

\section{A refined property of  the  interpolation operator $\Pi_h$}
This section   establishes a lower bound of $a_h(w-\Pi_hw,\Pi_hw)$. To this end,
 given any  element $K$,  we follow \cite{HuHuangLin10,HuShi12} to define $P_Kv\in P_4(K)$ by
\begin{equation}\label{QuasiInterpolation}
\int_{K}\na^{\ell}P_Kv dxdy=\int_K \na^{\ell}v dxdy, \ell=0, 1, 2, 3, 4,
\end{equation}
for any $v\in H^4(K)$.  Here and throughout this paper,  $\na^{\ell}v$ denotes the $\ell$-th order tensor of all
$\ell$-th order derivatives of $v$, for instance, $\ell=1$ the
gradient,  and $\ell=2$ the Hessian matrix, and that $\na_h^{\ell}$ are the
piecewise  counterparts of $\na^{\ell}$ defined element by element. Note that the operator $P_K$ is
well-posed.  It follows from the  definition of $P_K$ in
\eqref{QuasiInterpolation} that
\begin{equation}\label{commuting}
\na^4 P_K v=\Pi_{0,K}\na^4v
\end{equation}
with $\Pi_{0,K}$ the $L^2$ constant projection operator over $K$.
Then the global version $\Pi_0 $ of the interpolation operator
$\Pi_{0,K} $ is defined as
\begin{equation}
\Pi_0 |_K=\Pi_{0,K}  \text{ for any } K\in \cT_h.
\end{equation}

\begin{lemma}\label{Lemma2.1} For any $u\in P_4(K)$ and $v\in Q_{Ad}(K)$, there holds that
\begin{equation}\label{expansion}
\begin{split}
(\nabla^2(u-\Pi_Ku), \nabla^2 v)_{L^2(K)}&=
-\frac{h_{y, K}^2}{3}\int_K\frac{\pa^4u}{\pa x^2\pa y^2}\frac{\pa^2 v}{\pa x^2}dxdy\\
&\quad -\frac{h_{x, K}^2}{3}\int_K\frac{\pa^4u}{\pa x^2\pa
y^2}\frac{\pa^2 v}{\pa y^2}dxdy.
\end{split}
\end{equation}
\end{lemma}
\begin{proof}  Let $\xi$ and $\eta$ be defined as in \eqref{mapping}.
 It follows from the definition of $Q_{Ad}(K)$ that
\begin{equation}\label{eq2.10}
\begin{split}
&\frac{\pa^2 v}{\pa x^2}=a_0+a_1\xi+a_2\eta+a_3\xi\eta,\\
&\frac{\pa^2 v}{\pa y^2}=b_0+b_1\xi+b_2\eta+b_3\xi\eta,\\
&\frac{\pa^2 v}{\pa x\pa
y}=c_0+c_1\xi+c_2\eta+c_3\xi^2+c_4\eta^2,
\end{split}
\end{equation}
for some interpolation parameters $a_i$, $b_i$, $i=0, \cdots, 3$,  and $c_i$, $i=0,\cdots, 4$.
 Since $u\in P_4(K)$,  we have
 $$
 u=u_1+\frac{h_{x,K}^4}{4!}\frac{\pa^4 u}{\pa x^4}\xi^4
+\frac{h_{y,K}^4}{4!}\frac{\pa^4 u}{\pa y^4}\eta^4 
+\frac{h_{x,K}^2h_{y,K}^2}{4}\frac{\pa^4 u}{\pa x^2\pa
y^2}\xi^2\eta^2,
 $$
 where $u_1\in Q_{Ad}(K)$. Note that $\Pi_K u_1=u_1$,  and 
 $$
 \Pi_K\xi^4=2\xi^2-1,  \Pi_K\eta^4=2\eta^2-1, \text{ and } \Pi_K\xi^2\eta^2=\xi^2+\eta^2-1.
 $$ 
  This implies 
\begin{equation*}
\begin{split}
u-\Pi_Ku&=\frac{h_{x,K}^4}{4!}\frac{\pa^4 u}{\pa x^4}\big(\xi^2-1\big)^2
+\frac{h_{y,K}^4}{4!}\frac{\pa^4 u}{\pa y^4}\big(\eta^2-1\big)^2\\
&\quad +\frac{h_{x,K}^2h_{y,K}^2}{4}\frac{\pa^4 u}{\pa x^2\pa
y^2}\big(\xi^2-1\big) \big(\eta^2-1\big).
\end{split}
\end{equation*}
Therefore
\begin{equation}\label{eq2.11}
\begin{split}
&\frac{\pa^2 (u-\Pi_Ku)}{\pa x^2} =\frac{h_{x,K}^2}{4!}\frac{\pa^4
u}{\pa x^4}\big(12\xi^2-4\big)+\frac{h_{y,K}^2}{2}\frac{\pa^4 u}{\pa
x^2\pa
y^2}\big(\eta^2-1\big),\\
&\frac{\pa^2 (u-\Pi_Ku)}{\pa y^2} =\frac{h_{y,K}^2}{4!}\frac{\pa^4
u}{\pa y^4}\big(12\eta^2-4\big)+\frac{h_{x,K}^2}{2}\frac{\pa^4 u}{\pa
x^2\pa y^2}\big(\xi^2-1\big),\\
&\frac{\pa^2 (u-\Pi_Ku)}{\pa x\pa y} =h_{x,K}h_{y,K}\frac{\pa^4
u}{\pa x^2\pa y^2} \xi\eta.
\end{split}
\end{equation}
Since 
$$
\int_K \big(12\xi^2-4\big)(a_0+a_1\xi+a_2\eta+a_3\xi\eta)dxdy=0,$$
and
$$
\int_K\big(\eta^2-1\big)(a_1\xi+a_2\eta+a_3\xi\eta)dxdy=0,
$$
 a combination of \eqref{eq2.10} and \eqref{eq2.11} plus some
elementary calculation yield
\begin{equation}
\int_K\frac{\pa^2 (u-\Pi_Ku)}{\pa x^2}\frac{\pa^2 v}{\pa x^2}dxdy
=-\frac{h_{y,K}^2}{3}\int_K \frac{\pa^4 u}{\pa x^2\pa y^2}\frac{\pa^2
v}{\pa x^2}dxdy.
\end{equation}
A similar argument proves
\begin{equation}
\begin{split}
&\int_K\frac{\pa^2 (u-\Pi_Ku)}{\pa y^2}\frac{\pa^2 v}{\pa y^2}dxdy
=-\frac{h_{x,K}^2}{3}\int_K \frac{\pa^4 u}{\pa x^2\pa y^2}\frac{\pa^2
v}{\pa y^2}dxdy,\\
&\int_K\frac{\pa^2 (u-\Pi_Ku)}{\pa x\pa y }\frac{\pa^2 v}{\pa x \pa
y}dxdy =0,
\end{split}
\end{equation}
which completes the proof.
\end{proof}

%\begin{remark} After we finished the paper,  we found that a similar expansion of \eqref{expansion}
% was proposed and developed for several conforming/nonconforming  finite elements on rectangles
% for the eigenvalue problem of the second order Laplace  operator in a series of  works \cite{HuangLiLin2008,LinHuangLi2008,LinHuangLi09}.
%\end{remark}

The above lemma can be used to prove the following crucial lower bound.

\begin{lemma}\label{Lemma2.2} Suppose that $w\in H_0^2(\Omega)\cap H^4(\Omega)$ be the solution of  Problem \eqref{cont}.
Then,
\begin{equation}
\beta h^2 \leq (\nabla_h^2(w-\Pi_hw),
\nabla_h^2\Pi_hw)_{L^2(\Omega)},
\end{equation}
for some positive constant $\beta$ which is independent of the
meshsize $h$ provided that $\|f\|_{L^2(\Omega)}\not= 0$ and that the meshsize is small enough.
\end{lemma}
\begin{proof} Given $K\in\cT_h$, let the interpolation operator $P_K$
 be  defined as in \eqref{QuasiInterpolation}, which  leads to the
 following decomposition
 \begin{equation}\label{eq2.14}
 \begin{split}
  &(\nabla_h^2(w-\Pi_hw),
  \nabla_h^2\Pi_hw)_{L^2(\Omega)}\\
&=\sum\limits_{K\in\cT_h} (\nabla ^2(P_K w-\Pi_K P_Kw),
  \nabla^2\Pi_Kw)_{L^2(K)}\\
  &\quad +\sum\limits_{K\in\cT_h} (\nabla ^2(I-\Pi_K)(I-P_K)w,
  \nabla^2\Pi_Kw)_{L^2(K)}\\
&=I_1+I_2.
  \end{split}
 \end{equation}
 Let $u=P_Kw$ and $v=\Pi_kw$ in Lemma \ref{Lemma2.1}.  The first
 term $I_1$ on the right-hand side of \eqref{eq2.14} reads
\begin{equation*}
\begin{split}
I_1&=-\sum\limits_{K\in\cT_h}\frac{h_{y,K}^2}{3}\int_K\frac{\pa^4 P_Kw}{\pa x^2\pa y^2}\frac{\pa^2 \Pi_Kw}{\pa x^2}dxdy\\
&\quad -\sum\limits_{K\in\cT_h}\frac{h_{x,K}^2}{3}\int_K\frac{\pa^4
P_K w}{\pa x^2\pa y^2}\frac{\pa^2 \Pi_K w}{\pa y^2}dxdy,
\end{split}
\end{equation*}
which can be rewritten as
\begin{equation*}
\begin{split}
&I_1=-\sum\limits_{K\in\cT_h}\frac{h_{y,K}^2}{3}\int_K\frac{\pa^4 w}{\pa x^2\pa y^2}\frac{\pa^2 w}{\pa x^2}dxdy-\sum\limits_{K\in\cT_h}\frac{h_{x,K}^2}{3}\int_K\frac{\pa^4
 w}{\pa x^2\pa y^2}\frac{\pa^2  w}{\pa y^2}dxdy\\
&\quad +\sum\limits_{K\in\cT_h}\int_K\frac{\pa^4 (I-P_K)w}{\pa x^2\pa y^2}\bigg(\frac{h_{y,K}^2}{3}\frac{\pa^2 \Pi_Kw}{\pa x^2}+\frac{h_{x,K}^2}{3}\frac{\pa^2 \Pi_K w}{\pa y^2}\bigg)dxdy\\
&\quad+\sum\limits_{K\in\cT_h}\int_K\frac{\pa^4 w}{\pa x^2\pa y^2}\bigg(\frac{h_{y,K}^2}{3}\frac{\pa^2 (I-\Pi_K)w}{\pa x^2}+\frac{h_{x,K}^2}{3}\frac{\pa^2 (I-\Pi_K) w}{\pa y^2}\bigg)dxdy.
\end{split}
\end{equation*}
By the  commuting property of \eqref{commuting},
$$
\frac{\pa^4(I-P_K )w}{\pa x^2\pa y^2}=(I-\Pi_{0,K})\frac{\pa^4 w}{\pa x^2\pa y^2}.
$$
Note that 
$$
\|\frac{\pa^2 \Pi_K w}{\pa y^2}\|_{L^2(K)}+\|\frac{\pa^2 \Pi_K w}{\pa x^2}\|_{L^2(K)}\leq C \|w\|_{H^3{K}}.
$$
This, the error estimates of \eqref{AdinterpolationEstimate},   yield
\begin{equation*}
\begin{split}
I_1&=-\sum\limits_{K\in\cT_h}\frac{h_{y,K}^2}{3}\int_K\frac{\pa^4 w}{\pa x^2\pa y^2}\frac{\pa^2 w}{\pa x^2}dxdy\\
&\quad -\sum\limits_{K\in\cT_h}\frac{h_{x,K}^2}{3}\int_K\frac{\pa^4
w}{\pa x^2\pa y^2}\frac{\pa^2  w}{\pa
y^2}dxdy\\
&\quad +\mathcal{O}(h^2)\|(I-\Pi_0)\nabla^4 w\|_{L^2(\Omega)}\|w\|_{H^3(\Omega)}.
\end{split}
\end{equation*}
Since the mesh is uniform, an elementwise integration by parts yields
\begin{equation*}
\begin{split}
&-\sum\limits_{K\in\cT_h}\frac{h_{y,K}^2}{3}\int_K\frac{\pa^4 w}{\pa x^2\pa y^2}\frac{\pa^2 w}{\pa x^2}dxdy\\
&=\sum\limits_{K\in\cT_h}\frac{h_{y,K}^2}{3}\int_K\bigg(\frac{\pa^3 w}{\pa x^2\pa y}\bigg)^2 dxdy
-h_y^2\int_{\Gamma_y}\frac{\pa^3 w}{\pa x^2\pa y}\frac{\pa^2 w}{\pa x^2}\nu_2dx,
\end{split}
\end{equation*}
where $\Gamma_y$ is the boundary of $\Omega$ that parallels to the x-axis,  and $\nu_2$ is the second component of the unit normal vector 
$\nu=(\nu_1, \nu_2)$ of the boundary. Since $\frac{\pa  w}{ \pa y}=0$ on $\Gamma_y$,  $\frac{\pa^3 w}{\pa x^2\pa y}=0$ on $\Gamma_y$.
Hence,
\begin{equation*}
\begin{split}
-\sum\limits_{K\in\cT_h}\frac{h_{y,K}^2}{3}\int_K\frac{\pa^4 w}{\pa x^2\pa y^2}\frac{\pa^2 w}{\pa x^2}dxdy
=\sum\limits_{K\in\cT_h}\frac{h_{y,K}^2}{3}\int_K\bigg(\frac{\pa^3 w}{\pa x^2\pa y}\bigg)^2 dxdy
\end{split}
\end{equation*}
A similar procedure shows
$$
-\sum\limits_{K\in\cT_h}\frac{h_{x,K}^2}{3}\int_K\frac{\pa^4
w}{\pa x^2\pa y^2}\frac{\pa^2  w}{\pa
y^2}dxdy
=\sum\limits_{K\in\cT_h}\frac{h_{x,K}^2}{3}\int_K\bigg(\frac{\pa^3
w}{\pa x \pa y^2}\bigg)dxdy.
$$
Therefore
\begin{equation}\label{eq2.15}
\begin{split}
I_1&=\sum\limits_{K\in\cT_h}\frac{h_{y,K}^2}{3}\|\frac{\pa^3 w}{\pa x^2\pa
y}\|_{L^2(K)}^2 +\sum\limits_{K\in\cT_h} \frac{h_{x,K}^2}{3}\|\frac{\pa^3 w}{\pa x\pa
y^2}\|_{L^2(\Omega)}^2\\
&\quad +\mathcal{O}(h^2)\|(I-\Pi_0)\nabla^4 w\|_{L^2(\Omega)}\|w\|_{H^3(\Omega)}.
\end{split}
\end{equation}
The second term $I_2$ on the right-hand side of \eqref{eq2.14} can
be estimated by the error estimates of \eqref{AdinterpolationEstimate}
and the commuting property of \eqref{commuting}, which reads
\begin{equation}\label{eq2.16}
|I_2|\leq Ch^2 \|(I-\Pi_0)\nabla^4
w\|_{L^2(\Omega)}\|w\|_{H^3(\Omega)}.
\end{equation}
 Since the piecewise constant functions are dense
 in the space $L^2(\Om)$,
 $$
 \|(I-\Pi_0)\nabla^4
w\|_{L^2(\Omega)}\rightarrow 0 \text{ when } h\rightarrow 0.
 $$
 Since $\| f|_{L^2(\Omega)}\not =0$ implies that $|\frac{\pa^2 w}{\pa x \pa
y}|_{H^1(\Omega)}\not= 0$ (see more details in the following remark), a combination of
 \eqref{eq2.14}-\eqref{eq2.16} proves the desired result.
\end{proof}

\begin{remark}  For the rectangular domain $\Omega$ under consideration, the condition $|\frac{\pa^2 w}{\pa x \pa
y}|_{H^1(\Omega)}\not= 0$ holds provided that $\|f\|_{L^2(\Omega)}\not=0$.  In fact, if $|\frac{\pa^2 w}{\pa x \pa
y}|_{H^1(\Omega)}= 0$,  $w$ is of the form
$$
w=c_0xy+h(x)+g(y),
$$
for some function $h(x)$ with respect to $x$, and $g(y)$ with respect to $y$.  Then, the boundary condition concludes that
both $h(x)$ and $g(y)$ are constant. Hence the boundary condition indicates $w\equiv 0$, which contradicts with $w\not\equiv 0.$
\end{remark}

\begin{remark} The expansion \eqref{expansion} was analyzed in \cite{HuHuang11,LuoLin04,Yang2000}. Herein we give
 a new and much simpler proof. Moreover, compared with the
 regularity $H^5$ needed therein,  the analysis herein only needs the regularity $H^4$.
\end{remark}

\begin{remark} The  idea herein can be directly extended to the eigenvalue problem investigated in \cite{HuHuang11,Yang2000},
which improves and simplifies the analysis therein and  proves that the discrete eigenvalue produced
 by the Adini element is smaller  than the exact one  provided that the meshsize is sufficiently small.
  In addition, such a generalization  weakens the regularity condition from $u\in H^5(\Omega)$
 to $u\in H^4(\Omega)$  where $u$ is the eigenfunction. 
\end{remark}

\section{An identity of $(-f, w-w_h)$}
This section establishes the identity of $(-f, w-w_h)$ which is one man ingredient for the proof of Theorem \ref{main}. 

\begin{lemma}\label{Lemma3.1} Let $w$ and $w_h$ be solutions of problems \eqref{cont}
and \eqref{disc}, respectively. Then,
\begin{equation}\label{eq3.0}
\begin{split}
&(-f, w-w_h)_{L^2(\Omega)}\\
&=a_h(w,\Pi_hw-w_h)-(f, \Pi_hw-w_h)_{L^2(\Omega)}\\
&\quad +a_h(w-\Pi_hw,w-\Pi_hw)+a_h(w-\Pi_hw,w_h-\Pi_hw)\\
&\quad +2(f, \Pi_hw-w)_{L^2(\Omega)}+2a_h(w-\Pi_hw,\Pi_hw).
\end{split}
\end{equation}
\end{lemma}
\begin{proof}  We start with the following decomposition
\begin{equation}\label{eq3.1}
\begin{split}
&(-f, w-w_h)_{L^2(\Omega)}\\
&=(-f, w-w_h)+a_h(w,w-w_h)-a_h(w, w-w_h)\\
&=(-f,w-\Pi_hw)_{L^2(\Omega)}+(-f, \Pi_hw-w_h)_{L^2(\Omega)}\\
&\quad +a(w, \Pi_hw-w_h)+a_h(w, w-\Pi_hw)-a_h(w,w-w_h).
\end{split}
\end{equation}
The last two terms on the right-hand side of \eqref{eq3.1} allow for
a further decomposition:
\begin{equation}\label{eq3.2}
\begin{split}
&a_h(w, w-\Pi_hw)-a_h(w,w-w_h)\\
&=a_h(w-\Pi_hw, w-\Pi_hw)
+a_h(\Pi_hw,w-\Pi_hw)\\
&\quad -a_h(w-\Pi_hw,w-w_h)-a_h(\Pi_hw, w-w_h)\\
&=a_h(w-\Pi_hw, w_h-\Pi_hw) +a_h(\Pi_hw,w-\Pi_hw)\\
&\quad -a_h(\Pi_hw, w-w_h).
\end{split}
\end{equation}
 It follows from the discrete problem \eqref{disc} and the continuous problem
\eqref{cont} that the last term on the right-hand side of
\eqref{eq3.2} can be divided  as
\begin{equation}\label{eq3.3}
\begin{split}
&-a_h(\Pi_hw, w-w_h)\\
&=(f, \Pi_hw-w)_{L^2(\Omega)}-a_h(w, \Pi_hw-w)\\
&=(f, \Pi_hw-w)_{L^2(\Omega)}-a_h(w-\Pi_hw, \Pi_hw-w)\\
&\quad -a_h(\Pi_hw, \Pi_hw-w).
\end{split}
\end{equation}
A summary of \eqref{eq3.1}-\eqref{eq3.3} completes the proof.
\end{proof}

\begin{remark} The importance of the identity of \eqref{eq3.0} lies in that such a decomposition  separates the dominant term
$2a_h(w-\Pi_hw,\Pi_hw)$ from the other higher order terms,  which is the key  to employ  Lemma \ref{Lemma2.2}.
\end{remark}

\section{The conclusion and comments}
This paper presents the analysis of the $L^2$  norm error estimate of the Adini element.
It is proved that the best $L^2$ norm error estimate is at most of order $\mathcal{O}(h^2)$  which can not be
improved in general.  This result in fact indicates that the nonconforming Adini element space can not
contain any conforming space with an appropriate approximate property. This will cause further difficulty
for the a posteriori error analysis.    In fact, the reliable and efficient a posteriori error estimate for
 this element is still missing in the literature, see \cite{CGH12} for more details.

\end{document}